\documentclass[11pt,a4paper]{article}

\usepackage{inputenc}
\usepackage{amsmath}
\usepackage{bm}
\usepackage{bbold}
\usepackage{amsthm}

\usepackage{hyperref}

\title{Max-Plus Algebra Models of Queueing Networks\thanks{Proceedings of the International Workshop on Discrete Event Systems (WODES'96), University of Edinburgh, UK, August 19-21, 1996, London: IEE, 1996, pp.~76--81.}}

\author{Nikolai K. Krivulin\thanks{Faculty of Mathematics and Mechanics, St.~Petersburg State University, 28 Universitetsky Ave., St.~Petersburg, 198504, Russia, 
nkk@math.spbu.ru}
}

\date{}

\newtheorem{theorem}{Theorem}
\newtheorem{lemma}[theorem]{Lemma}

\def\sumo_#1^#2{\setbox0=\hbox{$\displaystyle{\sum}$}
                \setbox1=\hbox{$\scriptstyle{#1}$}
                \setbox2=\hbox{$\scriptstyle{#2}$}
		\setbox3=\hbox{${}_{{}_\oplus}\mathsurround=0pt$}
		\dimen1=.5\wd1 \advance\dimen1 by-.5\wd0
		\ifdim\dimen1>0pt
		   \ifdim\dimen1>\wd3 \kern\wd3 \else\kern\dimen1\fi\fi
		\dimen2=.5\wd2 \advance\dimen2 by-.5\wd0
		\ifdim\dimen2>0pt
		   \ifdim\dimen2>\wd3 \kern\wd3 \else\kern\dimen2\fi\fi
		\mathop{{\sum}{}_{{}_\oplus}}_{\kern-\wd3 #1}^{\kern-\wd3 #2}}
\def\sumol_#1{\setbox0=\hbox{$\displaystyle{\sum}$}
             \setbox1=\hbox{$\scriptstyle{#1}$}
	     \setbox3=\hbox{${}_{{}_\oplus}\mathsurround=0pt$}
	     \dimen1=.5\wd1 \advance\dimen1 by-.5\wd0
	     \ifdim\dimen1>0pt
	        \ifdim\dimen1>\wd3 \kern\wd3 \else\kern\dimen1\fi\fi
	     \mathop{{\sum}{}_{{}_\oplus}}_{\kern-\wd3 #1}}

\begin{document}

\maketitle

\begin{abstract}
A class of queueing networks which may have an arbitrary topology, and
consist of single-server fork-join nodes with both infinite and finite buffers
is examined to derive a representation of the network dynamics in terms of
max-plus algebra. For the networks, we present a common dynamic state equation
which relates the departure epochs of customers from the network nodes in an
explicit vector form determined by a state transition matrix. It is shown how
the matrices inherent in particular networks may be calculated from the
service times of customers. Since, in general, an explicit dynamic equation
may not exist for a network, related existence conditions are established in
terms of the network topology.
\\

\end{abstract}

\section{Introduction}

We consider a class of queueing networks with single-server nodes and
customers of a single class. In the networks, the server at each node is
supplied with a buffer which may have both infinite and finite capacity. There
is, in general, no restriction on the network topology; in particular, both
open and closed queueing networks may be included in the class.

In addition to the ordinary service procedure, specific fork-join operations
\cite{Bacc89,Gree91} may be performed in each node of the networks. In fact,
these operations allow customers (jobs, tasks) to be split into parts, and to
be merged into one, when circulating through the network. The fork-join
formalism proves to be useful in the description of dynamical processes in a
variety of actual systems, including production processes in manufacturing,
transmission of messages in communication networks, and parallel data
processing in multi-processor systems \cite{Bacc89}. As an example, one
can consider the splitting of a message into packets, and the merging of the
packets to restore the message, inherent in communication systems.

In this paper, the networks are examined so as to represent their dynamics in
terms of max-plus algebra \cite{Cuni79,Cohe89,Cass95}. The max-plus algebra
approach actually offers a quite compact and unified way of describing system
dynamics, which may provide a useful framework for analytical study and
computer simulation of discrete event systems including systems of queues.

It has been shown in \cite{Kriv94,Kriv95} that the evolution of both open and
closed tandem queueing systems may be described by the linear algebraic
equation
\begin{equation}\label{ve-sde}
\mbox{\boldmath $d$}(k) = T(k) \otimes \mbox{\boldmath $d$}(k-1),
\end{equation}
where $  \mbox{\boldmath $d$}(k)  $ is a vector of departure epochs from the
queues, $  T(k)  $ is a matrix calculated from service times of customers,
and $  \otimes  $ is an operator which determines the matrix-vector
multiplication in the max-plus algebra. In fact, this equation quite
frequently occurs in discrete event system analysis and simulation which are
based on the max-plus algebra approach. One can find a variety of related
examples in \cite{Cuni79,Cohe85,Bacc93,Cass95}.

The purpose of this paper is to show that the dynamics of the networks under
examination also allows of representation through dynamic state equation
(\ref{ve-sde}). We start with preliminary max-plus algebra definitions and
related results in Section~\ref{s-PDR}. Furthermore, Section~\ref{s-NM} gives
a general description of the network model, and shows how the dynamics of
nodes may be described through scalar equations in terms of max-plus algebra.
In Section~\ref{s-VR}, the scalar equations are extended to produce a vector
representation of the dynamics of the entire network. Finally, in
Section~\ref{s-ESE}, explicit dynamic state equations are derived in the form
of (\ref{ve-sde}). Since an explicit state dynamic equation does not have to
exist for an arbitrary network, related existence conditions in terms of
network topology are also included in Section~\ref{s-ESE}.

\section{Preliminary Definitions and Results}\label{s-PDR}

We start with a brief overview of basic algebraic facts and their graph
interpretation which we will exploit in the development of max-plus algebra
models of queueing networks. A detailed analysis of the max-plus algebra and
related algebraic systems, as well as their applications can be found in
\cite{Cuni79,Cohe85,Cohe89,Bacc93,Masl94,Cass95}.

{\it Max-plus algebra\/} is normally defined (see, e.g., \cite{Cass95}) as the
system $  \langle \underline{\mathbb{R}}, \oplus, \otimes \rangle $, where
$  \underline{\mathbb{R}} = \mathbb{R} \cup \{\varepsilon\}  $ is the
set of real numbers with $  \varepsilon = -\infty  $ adjoined, and the
symbols $  \oplus  $ and $  \otimes  $ present binary operations
determined for any $  x,y \in \underline{\mathbb{R}}  $ respectively as
$$
x \oplus y = \max(x,y), \quad x \otimes y = x + y.
$$

As one can verify \cite{Cuni79,Cohe89}, most of the properties of the ordinary
addition and multiplication, including their associativity and commutativity,
as well as distributivity of multiplication over addition, are extended to the
operations $  \oplus  $ and $  \otimes $. These properties allow usual
algebraic manipulations in the max-plus algebra to be performed under the
standard conventions regarding brackets and precedence of multiplication over
addition. Note that, in contrast to the conventional algebra, the operation
$  \oplus  $ is idempotent; that is, for any
$  x \in \underline{\mathbb{R}} $, we have $  x \oplus x = x $.

There are the null and identity elements in the max-plus algebra, namely
$  \varepsilon  $ and $  e = 0 $, to satisfy the conditions
$  x \oplus \varepsilon = \varepsilon \oplus x = x $, and
$  x \otimes e = e \otimes x = x $, for any
$  x \in \underline{\mathbb{R}} $. The absorption rule which involves
$  x \otimes \varepsilon = \varepsilon \otimes x = \varepsilon  $ is also
true in this algebra.

The max-plus algebra of matrices is introduced in the regular way
\cite{Cuni79,Cohe85,Cohe89}. Specifically, for any $(n \times n)$-matrices
$  X = (x_{ij})  $ and $  Y = (y_{ij}) $, the entries of
$  U = X \oplus Y  $ and $  V = X \otimes Y  $ are calculated as
$$
u_{ij} = x_{ij} \oplus y_{ij}, \quad \mbox{and} \quad
v_{ij} = \sumo_{k=1}^{n} x_{ik} \otimes y_{kj},
$$
where the symbol $  \sum_{{}^{\oplus}}  $ denotes the iterated operation
$  \oplus $. As the null element, the matrix $  \mathcal{E}  $ with all its
entries equal to $  \varepsilon  $ is taken in this algebra, while the
matrix $  E = {\rm diag}(e,\ldots,e)  $ with the off-diagonal entries equal
to $  \varepsilon  $ presents the identity element.

In perfect analogy to the conventional matrix algebra, one can define for any
square matrix $  X $,
$$
X^{0} = E, \quad
X^{q}
= \underbrace{X \otimes \cdots \otimes X}_{\mbox{\scriptsize $ q  $ times}}
\quad \mbox{for $  q=1,2,\ldots $}
$$
Note, however, that idempotency in the max-plus algebra leads, in particular,
to the identity \cite{Cuni79}
$$
(E \oplus X)^{q} = E \oplus X \oplus \cdots \oplus X^{q}.
$$

Many phenomena inherent in the matrix max-plus algebra appear to be well
explained in terms of their graph interpretations
\cite{Cuni79,Cohe85,Masl94,Cass95}. To illustrate, we can consider an
$(n \times n)$-matrix $  X  $ with its entries
$  x_{ij} \in \underline{\mathbb{R}} $, and note that it can be treated as
the adjacency matrix of an oriented graph with $  n  $ nodes, provided each
entry $  x_{ij} \neq \varepsilon  $ implies the existence of the arc
$  (i,j)  $ in the graph, whereas $  x_{ij} = \varepsilon  $ does the lack
of the arc. The graph is then said to be {\it associated with\/} the matrix
$  X $.

Let us calculate the matrix $  X^{2} = X \otimes X $, and denote its entries
by $  x_{ij}^{(2)} $. Clearly, we have $  x_{ij}^{(2)} \neq \varepsilon  $
if and only if there exists at least one path from node $  i  $ to node
$  j  $ in the graph, which consists of two arcs. Moreover, for any positive
integer $  q $, the matrix $  X^{q}  $ has the entry
$  x_{ij}^{(q)} \neq \varepsilon  $ only when there exists a path with the
length $  q  $ from $  i  $ to $  j $.

Suppose that the graph associated with the matrix $  X  $ is acyclic. It is
clear that we will have $  X^{q} = \mathcal{E}  $ for all $  q > p $, where
$  p  $ is the length of the longest path in the graph. Assume now the graph
not to be acyclic, and then consider any one of its circuits. Since it is
possible to construct a cyclic path of any length, which lies along the
circuit, we conclude that $  X^{q} \neq \mathcal{E}  $ for all
$  q=1,2,\ldots $

Finally, we consider the implicit equation in the unknown vector
$  \mbox{\boldmath $x$} = (x_{1},\ldots,x_{n})^{T} $,
\begin{equation}\label{e-I}
\mbox{\boldmath $x$} = U \otimes \mbox{\boldmath $x$}
                     \oplus \mbox{\boldmath $v$},
\end{equation}
where $  U = (u_{ij})  $ and
$  \mbox{\boldmath $v$} = (v_{1},\ldots,v_{n})^{T}  $ are respectively given
$(n\times n)$-matrix and $n$-vector. This equation actually plays a large
role in max-plus algebra representations of dynamical systems including
systems of queues \cite{Cohe85,Cass95,Kriv95}. The next lemma offers
particular conditions for (\ref{e-I}) to be solvable, and shows how the
solution may be calculated. One can find a detailed analysis of (\ref{e-I}) in
the general case in \cite{Cohe89}.

\begin{lemma}\label{l-1}
Suppose that the entries of the matrix $  U  $ and the vector
$  \mbox{\boldmath $v$}  $ are either positive or equal to
$  \varepsilon $. Then equation {\rm (\ref{e-I})} has the unique bounded
solution $  \mbox{\boldmath $x$}  $ if and only if the graph associated with
$  U  $ is acyclic. Provided that the solution exists, it is given by
$$
\mbox{\boldmath $x$} = (E \oplus U)^{p} \otimes \mbox{\boldmath $v$},
$$
where $  p  $ is the length of the longest path in the graph.
\end{lemma}

To prove the lemma, first note that recurrent substitution of
$  \mbox{\boldmath $x$}  $ from equation (\ref{e-I}) into its right-hand
side, made $  q  $ times, and trivial manipulations give
$$
\mbox{\boldmath $x$} = U^{q+1} \otimes \mbox{\boldmath $x$}
\oplus (E \oplus U \oplus \cdots \oplus U^{q}) \otimes \mbox{\boldmath $v$}.
$$

The rest of the proof may be readily furnished based on the above graph
interpretation as well as on the idempotency of the operation $  \oplus $.

\section{The Network Model}\label{s-NM}

In this section, we present a network model which may be considered as an
extension of acyclic fork-join queueing networks investigated in
\cite{Bacc89,Gree91}. In fact, we do not restrict ourselves on acyclic
networks, but assume the networks to have an arbitrary topology. Moreover, we
examine not only the networks with the infinite capacity of buffers in their
nodes, but also those with finite buffers and blocking of servers.

\subsection{A General Description of the Model}

\begin{figure}[hhh]
\begin{center}
\begin{picture}(80,35)

\put(5.5,18){$1$}
\put(0,9){$r_{1}\!=\infty$}
\put(0,5.5){$s_{1}\!=\infty$}
\put(0,16){\thicklines\line(1,0){4}}
\put(0,12){\thicklines\line(1,0){4}}
\put(4,16){\thicklines\line(0,-1){4}}
\put(6,14){\thicklines\circle{3}}

\put(8,14){\vector(1,0){9}}

\put(22.5,18){$2$}
\put(16,9){$r_{2},s_{2}$}
\put(17,16){\thicklines\line(1,0){4}}
\put(17,12){\thicklines\line(1,0){4}}
\put(21,16){\thicklines\line(0,-1){4}}
\put(23,14){\thicklines\circle{3}}

\put(25,14){\vector(1,1){9}}
\put(25,14){\vector(1,-1){9}}

\put(39.5,27){$3$}
\put(33,18){$r_{3},s_{3}$}
\put(34,25){\thicklines\line(1,0){4}}
\put(34,21){\thicklines\line(1,0){4}}
\put(38,25){\thicklines\line(0,-1){4}}
\put(40,23){\thicklines\circle{3}}

\put(42,23){\vector(1,0){18}}
\put(42,23){\vector(1,-1){18}}

\put(65.5,27){$4$}
\put(59,18){$r_{4},s_{4}$}
\put(60,25){\thicklines\line(1,0){4}}
\put(60,21){\thicklines\line(1,0){4}}
\put(64,25){\thicklines\line(0,-1){4}}
\put(66,23){\thicklines\circle{3}}

\put(68,23){\line(1,0){9}}
\put(77,23){\line(0,1){10}}
\put(77,33){\line(-1,0){67}}
\put(10,33){\line(0,-1){12}}
\put(10,21){\vector(1,-1){7}}

\put(39.5,9){$5$}
\put(33,0){$r_{5},s_{5}$}
\put(34,7){\thicklines\line(1,0){4}}
\put(34,3){\thicklines\line(1,0){4}}
\put(38,7){\thicklines\line(0,-1){4}}
\put(40,5){\thicklines\circle{3}}

\put(42,5){\vector(1,0){18}}

\put(65.5,9){$6$}
\put(59,0){$r_{6},s_{6}$}
\put(60,7){\thicklines\line(1,0){4}}
\put(60,3){\thicklines\line(1,0){4}}
\put(64,7){\thicklines\line(0,-1){4}}
\put(66,5){\thicklines\circle{3}}

\put(68,5){\vector(1,0){9}}

\end{picture}

\end{center}
\caption{A queueing network with $  n=6  $ nodes.}\label{f-FJQN}

\end{figure}
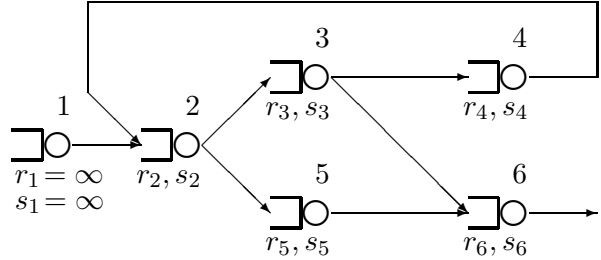

We consider a queueing network consisting of $  n  $ single-server nodes,
with customers of a single class, which circulate through the network. An
example of the network under study is shown in Fig.~\ref{f-FJQN}.

The topology of the network is described by an oriented graph
$  \mathcal{G} = ({\bf N}, {\bf A})  $ which, in general, does not have to be
acyclic. In the graph $  \mathcal{G} $, the set $  {\bf N} = \{1,\ldots,n\}  $
represents the nodes of the network, and the set
$  {\bf A} = \{(i,j)\} \subset {\bf N} \times {\bf N}  $ consists of arcs
determining the transition routes of customers. For any
$  i,j \in {\bf N} $, the arc $  (i,j)  $ belongs to $  {\bf A}  $ if
and only if the $i$th node passes customers directly to node $  j $.

For every node $  i \in {\bf N} $, we introduce the set of its predecessors
$  {\bf P}(i) = \{ j | \, (j,i) \in {\bf A} \}  $ and the set of its
successors $  {\bf S}(i) = \{ j | \, (i,j) \in {\bf A} \} $. We suppose that,
in specific cases, there may be one of the conditions
$  {\bf P}(i) = \emptyset  $ and $  {\bf S}(i) = \emptyset  $ encountered.
Each node $  i  $ with $  {\bf P}(i) = \emptyset  $ is assumed to
represent an infinite external arrival stream of customers. Provided that
$  {\bf S}(i) = \emptyset $, the node is considered as an output node
intended to release customers from the network.

Each node $  i \in {\bf N}  $ includes a server and its buffer which
together present a single-server queue operating under the first-come,
first-served (FCFS) queueing discipline. The buffer at the server in node
$  i  $ may have either finite or infinite capacity $  s_{i} $; that is,
$  0 \leq s_{i} \leq \infty $. At the initial time, the server at each node
$  i  $ is assumed to be free of customers, whereas in its buffer, there may
be $  r_{i} $, $  0 \leq r_{i} \leq s_{i} $, customers waiting for service.
It is thought that the values $  s_{i} = r_{i} = \infty  $ are set for every
node $  i  $ with $  {\bf P}(i) = \emptyset $, representing an external
arrival stream of customers. We consider the numbers $  r_{i}  $ and
$  s_{i} $, $  i=1,\ldots,n $, as initial conditions in the model.

To describe the dynamics of the queue in node $  i $, we use the following
symbols:
\begin{list}{}{\topsep=0pt\itemsep=0pt\itemindent\parindent}
\item[$ a_{i}(k), $] the $k$th arrival epoch to the queue;
\item[$ b_{i}(k), $] the $k$th service initiation time;
\item[$ c_{i}(k), $] the $k$th service completion time;
\item[$ d_{i}(k), $] the $k$th departure epoch from the queue.
\end{list}
Furthermore, the service time of the $k$th customer at server $  i  $ is
denoted by $  \tau_{ik} $, $  \tau_{ik} > 0 $. We assume that
$  \tau_{ik}  $ are given parameters for all $  i=1,\dots,n $, and
$  k=1,2,\dots $, while $  a_{i}(k) $, $  b_{i}(k) $, $  c_{i}(k) $, and
$  d_{i}(k)  $ are considered as unknown state variables. In addition, with
the condition that the network starts operating at time zero, it is convenient
to set $  d_{i}(0) \equiv e $, and $  d_{i}(k) \equiv \varepsilon  $ for
all $  k < 0 $, $  i=1,\ldots,n $.

As one can see, relations between the state variables, which are actually
determined by the network topology, initial conditions, and special features
inherent in node operation, just represent the dynamics of the network. We
will describe the network dynamics in more detail and give related algebraic
representations in the subsequent sections.

\subsection{The Dynamics of Nodes}

We suppose that, in addition to the usual service procedure, special join and
fork operations are performed in the nodes, respectively before and after
service \cite{Bacc89}. With the condition that all buffers at servers in a
network have infinite capacity, the fork-join mechanism may be described as
follows.

The {\it join\/} operation is actually thought to cause each customer which
comes into node $  i $, not to enter the buffer at the server but to wait
until at least one customer from every node $  j \in {\bf P}(i)  $ arrives.
As soon as these customers arrive, they, taken one from each preceding node,
are united to be treated as being one customer which then enters the buffer to
become a new member of the queue. Note that only the customers who are waiting
for service may be placed into the buffer at the node. Those customers which
are ready to be joined, but, in the absence of all required customers, have
not been joined yet, are assumed to reside in another place, say in an
auxiliary buffer available at the node. It is suggested that the auxiliary
buffers invariably have infinite capacity.

The {\it fork\/} operation at node $  i  $ is initiated every time a
customer releases the server after completion of his service; it consists in
giving rise to several new customers to substitute for the former one. As many
new customers appear in node $  i  $ as there are succeeding nodes in the
set $  {\bf S}(i) $. The customers simultaneously depart the node, each going
to separate node $  j \in {\bf S}(i) $. Finally, we assume that the execution
of the fork-join operations when appropriate customers are available, as well
as the transition of customers within and between nodes require no time.

It is easy to see that, assuming the sets $  {\bf P}(i)  $ and
$  {\bf S}(i)  $ to include no more than one node for each
$  i \in {\bf N} $, one can arrive at a queueing system in which essentially
no fork-join operations are performed. As examples, both open and closed
tandem systems may be considered \cite{Kriv94,Kriv95}, which actually present
queueing networks with the simplest topology.

In order to set up the equations which represent the dynamics of nodes, let us
first consider a network with infinite buffers. It follows from the above
description of the fork-join mechanism that the time of the $k$th arrival into
the queue at node $  i $, which actually coincides with that of the
completion of the $k$th join operation, may be represented as
\cite{Bacc89,Gree91}
\begin{equation}\label{e-a1}
a_{i}(k) = \left\{
	     \begin{array}{ll}
	       \displaystyle{\sumol_{j\in \mbox{\scriptsize\bf P}(i)}}
                                                               d_{j}(k-r_{i}),
                               & \mbox{if $  {\bf P}(i) \neq \emptyset $}, \\
               \quad\varepsilon, & \mbox{if $  {\bf P}(i) = \emptyset $},
	     \end{array}
	   \right.
\end{equation}
whereas the equations which determine the other state variables are readily
written in the form
\begin{eqnarray}
b_{i}(k) & = & a_{i}(k) \oplus d_{i}(k-1), \label{e-b1} \\
c_{i}(k) & = & \tau_{ik} \otimes b_{i}(k), \label{e-c1} \\
d_{i}(k) & = & c_{i}(k). \label{e-d1}
\end{eqnarray}

Suppose now that the buffers at servers in the network may have limited
capacity. In such systems, servers may be blocked according to some blocking
mechanism \cite{Gree91,Chen94}. The rest of the section shows how to
represent the dynamics of nodes operating under the {\it manufacturing\/} and
{\it communication\/} blocking rules, both being commonly encountered in
practice.

Let us first assume the network operation to follow the manufacturing blocking
rule. Application of this type of blocking implies that, upon completion of
his service at node $  i $, a customer cannot release the server at the node
if there is at least one succeeding node $  j \in {\bf S}(i)  $ in which the
buffer is full. As soon as all nodes included in $  {\bf S}(i)  $ regain an
empty buffer space, the customer leaves the server to produce new customers
which have to depart node $  i  $ immediately.

The inclusion of manufacturing blocking leads us to the new equation
representing departure times, which is to substitute for (\ref{e-d1}),
\begin{equation}\label{e-d2}
d_{i}(k) = c_{i}(k) \oplus {\cal D}_{i}(k),
\end{equation}
where
$$
{\cal D}_{i}(k)
= \left\{
    \begin{array}{ll}
      \displaystyle{\sumol_{j \in \mbox{\scriptsize\bf S}(i)}}
                                                             d_{j}(k-s_{j}-1),
                               & \mbox{if $  {\bf S}(i) \neq \emptyset $}, \\
      \quad\varepsilon,        & \mbox{if $  {\bf S}(i) = \emptyset $}.
    \end{array}
   \right.
$$
Clearly, equations (\ref{e-a1}--\ref{e-c1}) remain unchanged.

Finally, we suppose that the network operates under communication blocking.
This blocking rule requires the server in node $  i  $ not to initiate
service of a customer until there is an empty space in the buffer in each node
$  j \in {\bf S}(i) $. To represent the dynamics of node $  i $, one may
take equations (\ref{e-a1}), (\ref{e-c1}), and (\ref{e-d1}) respectively for
$  a_{i}(k) $, $  c_{i}(k) $, and $  d_{i}(k) $. With the symbol
$  {\cal D}_{i}(k)  $ introduced above, an appropriate equation for
$  b_{i}(k)  $ is now written as
\begin{equation}\label{e-b2}
b_{i}(k) = a_{i}(k) \oplus d_{i}(k-1) \oplus {\cal D}_{i}(k).
\end{equation}

\section{A Vector Representation}\label{s-VR}

We now turn to the algebraic representation of the dynamics of the entire
network. To describe the dynamics in a compact form, we introduce the vectors
\begin{eqnarray*}
\mbox{\boldmath $a$}(k) & = & (a_{1}(k),\ldots,a_{n}(k))^{T}, \\
\mbox{\boldmath $b$}(k) & = & (b_{1}(k),\ldots,b_{n}(k))^{T}, \\
\mbox{\boldmath $c$}(k) & = & (c_{1}(k),\ldots,c_{n}(k))^{T}, \\
\mbox{\boldmath $d$}(k) & = & (d_{1}(k),\ldots,d_{n}(k))^{T},
\end{eqnarray*}
and the diagonal matrix
$$
{\cal T}_{k}
= \left(
    \begin{array}{ccc}
      \tau_{1k}   &        & \varepsilon \\
                  & \ddots &             \\
      \varepsilon &        & \tau_{nk}
    \end{array}
  \right).
$$

\subsection{Networks with Infinite Buffers}\label{s-NIB}

We start with the derivation of a vector representation relevant to equations
(\ref{e-a1}--\ref{e-d1}) set up for networks with infinite buffers. First note
that vector equations associated with (\ref{e-b1}--\ref{e-d1}) may be written
immediately.

To get equation (\ref{e-a1}) in a vector form, we define
$  M_{r} = \max \{ r_{i} | \, r_{i} < \infty, \, i=1,\ldots,n \} $. It is
easy to see that we may now represent (\ref{e-a1}) as
$$
a_{i}(k) = \sumo_{m=0}^{M_{\scriptstyle r}}
                                \sumo_{j=1}^{n} g_{ji}^{m} \otimes d_{j}(k-m),
$$
where the numbers $  g_{ij}^{m}  $ are determined using the topology of the
network by the condition
$$
g_{ij}^{m} = \left\{
	       \begin{array}{ll}
                e, & \mbox{if $  i \in {\bf P}(j)  $ and $  m=r_{j} $}, \\
                \varepsilon, & \mbox{otherwise}.
               \end{array}
	     \right.
$$
Furthermore, we introduce the matrices $  G_{m} = \left(g_{ij}^{m}\right)  $
for each $  m=0,1,\ldots,M_{r} $, and then bring the above equation into its
associated vector form
$$
\mbox{\boldmath $a$}(k)
= \sumo_{m=0}^{M_{\scriptstyle r}} G_{m}^{T}
                                            \otimes \mbox{\boldmath $d$}(k-m),
$$
where $  G_{m}^{T}  $ denotes the transpose of the matrix $  G_{m} $. Note
that each matrix $  G_{m}  $ presents an adjacency matrix of the partial
graph $  \mathcal{G}_{m} = ({\bf N},{\bf A}_{m})  $ with
$  {\bf A}_{m} = \{(i,j) | \, i \in {\bf P}(j), \, r_{j}=m\} $.

We are now in a position to describe the network dynamics in vector terms. By
replacing equations (\ref{e-a1}--\ref{e-d1}) with their vector
representations, we obtain
\begin{eqnarray*}
\mbox{\boldmath $a$}(k) & = & \sumo_{m=0}^{M_{\scriptstyle r}} G_{m}^{T}
                                         \otimes \mbox{\boldmath $d$}(k-m), \\
\mbox{\boldmath $b$}(k) & = & \mbox{\boldmath $a$}(k)
                                          \oplus \mbox{\boldmath $d$}(k-1), \\
\mbox{\boldmath $c$}(k) & = & {\cal T}_{k} \otimes \mbox{\boldmath $b$}(k), \\
\mbox{\boldmath $d$}(k) & = & \mbox{\boldmath $c$}(k).
\end{eqnarray*}
Clearly, these equations can be reduced to an equation in one vector variable,
say $  \mbox{\boldmath $d$}(k) $. In that case, appropriate substitutions
will lead us to the equation
\begin{eqnarray}
\mbox{\boldmath $d$}(k) & = & {\cal T}_{k} \otimes G_{0}^{T}
                                            \otimes \mbox{\boldmath $d$}(k)
            \oplus {\cal T}_{k} \otimes \mbox{\boldmath $d$}(k-1) \nonumber \\
&   & \mbox{} \oplus {\cal T}_{k}
                            \otimes \sumo_{m=1}^{M_{\scriptstyle r}} G_{m}^{T}
                               \otimes \mbox{\boldmath $d$}(k-m). \label{vi-i}
\end{eqnarray}

\subsection{Networks with Finite Buffers}

Consider a network with finite buffers, and assume that it operates under
the manufacturing blocking rule. With
$  M_{s} = \max \{ s_{i}+1 | \, s_{i} < \infty, \, i=1,\ldots,n \} $,
equation (\ref{e-d2}) may be put in the form
$$
d_{i}(k) = c_{i}(k) \oplus \sumo_{m=1}^{M_{\scriptstyle s}}
                                \sumo_{j=1}^{n} h_{ij}^{m} \otimes d_{j}(k-m),
$$
where
\begin{equation}\label{e-h}
h_{ij}^{m} = \left\{
	       \begin{array}{ll}
                e, & \mbox{if $  j \in {\bf S}(i)  $ and $  m=s_{j}+1 $},\\
                \varepsilon, & \mbox{otherwise}.
               \end{array}
	     \right.
\end{equation}
In a similar way as in Section~\ref{s-NIB}, one can introduce the matrices
$  H_{m} = \left(h_{ij}^{m}\right) $, $  m=1,\ldots, M_{s} $, and then
rewrite (\ref{e-d2}) so as to get a representation for
$  \mbox{\boldmath $d$}(k) $.

Taking into account that the equations set up previously to represent the
vectors $  \mbox{\boldmath $a$}(k) $, $  \mbox{\boldmath $b$}(k) $, and
$  \mbox{\boldmath $c$}(k)  $ remain valid, we arrive at the set of
equations
\begin{eqnarray*}
\mbox{\boldmath $a$}(k) & = & \sumo_{m=0}^{M_{\scriptstyle r}} G_{m}^{T}
                                         \otimes \mbox{\boldmath $d$}(k-m), \\
\mbox{\boldmath $b$}(k) & = & \mbox{\boldmath $a$}(k)
                                          \oplus \mbox{\boldmath $d$}(k-1), \\
\mbox{\boldmath $c$}(k) & = & {\cal T}_{k}
                                           \otimes \mbox{\boldmath $b$}(k), \\
\mbox{\boldmath $d$}(k) & = & \mbox{\boldmath $c$}(k)
                              \oplus \sumo_{m=1}^{M_{\scriptstyle s}} H_{m}
                                            \otimes \mbox{\boldmath $d$}(k-m).
\end{eqnarray*}

Without loss of generality, we consider that $  M_{r}=M_{s}=M $. If it
actually holds that $  M_{r} < M_{s} $ (the inequality $  M_{r} > M_{s}  $
is contradictory to the initial conditions), one may set $  M = M_{s} $, and
then define $  G_{m} = \mathcal{E}  $ for all
$  m=M_{r}+1,M_{r}+2,\ldots,M_{s} $. With this assumption, we may drop
the subscripts so as to write $  M  $ instead of both $  M_{r}  $ and
$  M_{s} $.

Proceeding to an equation in $  \mbox{\boldmath $d$}(k) $, we get
\begin{eqnarray}
\lefteqn{\mbox{\boldmath $d$}(k) = {\cal T}_{k} \otimes G_{0}^{T}
                                            \otimes \mbox{\boldmath $d$}(k)
           \oplus {\cal T}_{k} \otimes \mbox{\boldmath $d$}(k-1)} \nonumber \\
& & \mbox{} \oplus \sumo_{m=1}^{M} ({\cal T}_{k} \otimes G_{m}^{T}
            \oplus H_{m}) \otimes \mbox{\boldmath $d$}(k-m). \label{vi-m}
\end{eqnarray}

Let us now assume the network to follow the communication blocking rule. In
the same way as for manufacturing blocking, one may define matrices
$  H_{1},\ldots,H_{M}  $ through (\ref{e-h}), and represent equation
(\ref{e-b2}) in its vector form. The set of vector equations describing the
network dynamics then becomes
\begin{eqnarray*}
\mbox{\boldmath $a$}(k) & = & \sumo_{m=0}^{M} G_{m}^{T}
                                         \otimes \mbox{\boldmath $d$}(k-m), \\
\mbox{\boldmath $b$}(k) & = & \mbox{\boldmath $a$}(k)
                                              \oplus \mbox{\boldmath $d$}(k-1)
            \oplus \sumo_{m=1}^{M} H_{m} \otimes \mbox{\boldmath $d$}(k-m), \\
\mbox{\boldmath $c$}(k) & = & {\cal T}_{k} \otimes \mbox{\boldmath $b$}(k), \\
\mbox{\boldmath $d$}(k) & = & \mbox{\boldmath $c$}(k).
\end{eqnarray*}

Finally, by combining these equations, we have
\begin{eqnarray}
\lefteqn{\mbox{\boldmath $d$}(k) = {\cal T}_{k} \otimes G_{0}^{T}
                                            \otimes \mbox{\boldmath $d$}(k)
           \oplus {\cal T}_{k} \otimes \mbox{\boldmath $d$}(k-1)} \nonumber \\
& & \mbox{} \oplus {\cal T}_{k} \otimes \sumo_{m=1}^{M} (G_{m}^{T}
            \oplus H_{m}) \otimes \mbox{\boldmath $d$}(k-m). \label{vi-c}
\end{eqnarray}

\section{The Explicit State Equation}\label{s-ESE}

Let us consider equations (\ref{vi-i}), (\ref{vi-m}), and (\ref{vi-c}) derived
above, and note that they actually present implicit equations in the system
state variable $  \mbox{\boldmath $d$}(k) $. In this section, we show how
these equations may be put in their associated explicit forms which are
normally more suitable for analytical treatments and computer simulation of
the network dynamics. Since, in general, the implicit equations do not have to
be explicitly solvable, the conditions for an explicit state equation to exist
are also established.

In order to examine the implicit equations, we first note that they all take
the form of (\ref{e-I}) with $  U = {\cal T}_{k} \otimes G_{0}^{T} $. Since
the matrix $  {\cal T}_{k}  $ is diagonal, each graph associated with the
matrix $  G_{0}^{T}  $ will be likewise associated with $  U $. In
addition, the graph associated with the matrix $  G_{0}  $ and that with its
transpose are both acyclic or not at once. Finally, both graphs have a common
length $  p  $ of their longest paths.

Now it is not difficult to apply Lemma~\ref{l-1} so as to prove the following
statement.
\begin{theorem}\label{t-1}
Suppose that in the network model with infinite buffers, the graph
$  \mathcal{G}_{0}  $ associated with the matrix $  G_{0}  $ is acyclic. Then
equation {\rm (\ref{vi-i})} can be solved to produce the explicit state
dynamic equation
\begin{equation}\label{ve-de}
\mbox{\boldmath $d$}(k) = \sumo_{m=1}^{M} T_{m}(k)
                                            \otimes \mbox{\boldmath $d$}(k-m),
\end{equation}
with the state transition matrices
\begin{eqnarray*}
T_{1}(k) & = & (E \oplus {\cal T}_{k} \otimes G_{0}^{T})^{p}
                         \otimes {\cal T}_{k} \otimes (E \oplus G_{1}^{T}), \\
T_{m}(k) & = & (E \oplus {\cal T}_{k} \otimes G_{0}^{T})^{p}
                                    \otimes {\cal T}_{k} \otimes G_{m}^{T}, \\
         &   & \mbox{} m=2,\ldots,M,
\end{eqnarray*}
where $  p  $ is the length of the longest path in $  \mathcal{G}_{0} $.
\end{theorem}

As one can see, the matrix coefficient at $  \mbox{\boldmath $d$}(k)  $ on
the right-hand side of both equation (\ref{vi-m}) and (\ref{vi-c}), which is
just responsible for explicit representation (\ref{ve-de}) to exist, remains
the same as in equation (\ref{vi-i}) set up for the network with infinite
buffers. We may therefore conclude that entering finite buffers into the
network model has no effect on the existence of its associated explicit
dynamic state equation.

Now one can reformulate Theorem~\ref{t-1} to extend it to the networks with
finite buffers. In short, under the same conditions as in the theorem,
equations (\ref{vi-m}) and (\ref{vi-c}) may be put in the form of
(\ref{ve-de}), with the state transition matrices defined respectively as
\begin{eqnarray*}
T_{1}(k) \!\! & = & \!\! (E \oplus {\cal T}_{k} \otimes G_{0}^{T})^{p}
 \otimes ({\cal T}_{k} \oplus {\cal T}_{k} \otimes G_{1}^{T} \oplus H_{1}), \\
T_{m}(k) \!\! & = & \!\! (E \oplus {\cal T}_{k} \otimes G_{0}^{T})^{p}
                     \otimes ({\cal T}_{k} \otimes G_{m}^{T} \oplus H_{m}), \\
         &   & \mbox{} m=2,\ldots,M,
\end{eqnarray*}
and
\begin{eqnarray*}
T_{1}(k) \! & = & \! (E \oplus {\cal T}_{k} \otimes G_{0}^{T})^{p}
            \otimes {\cal T}_{k} \otimes (E \oplus G_{1}^{T} \oplus H_{1}), \\
T_{m}(k) \! & = & \! (E \oplus {\cal T}_{k} \otimes G_{0}^{T})^{p}
            \otimes {\cal T}_{k} \otimes (G_{m}^{T} \oplus H_{m}), \\
         &   & \mbox{} m=2,\ldots,M.
\end{eqnarray*}

Finally, with the extended state vector
$$
\widehat{\mbox{\boldmath $d$}}(k) = \left(
                                      \begin{array}{l}
                                        \mbox{\boldmath $d$}(k) \\
                                        \mbox{\boldmath $d$}(k-1) \\
                                        \vdots \\
                                        \mbox{\boldmath $d$}(k-M+1)
                                      \end{array}
                                    \right),
$$
we may bring (\ref{ve-de}) into the form of (\ref{ve-sde}):
$$
\widehat{\mbox{\boldmath $d$}}(k)
                 = \widehat{T}(k) \otimes \widehat{\mbox{\boldmath $d$}}(k-1),
$$
where the new state transition matrix is defined as
$$
\widehat{T}(k) = \left(
		   \begin{array}{ccccc}
                     T_{1}(k) & T_{2}(k) & \cdots & \cdots &T_{M}(k)  \\
		     E        & \mathcal{E} & \cdots & \cdots & \mathcal{E} \\
		              & \ddots   & \ddots &        & \vdots   \\
		              &          & \ddots & \ddots & \vdots   \\
		     \mathcal{E} &          &        & E      & \mathcal{E}
		   \end{array}
		 \right).
$$

In conclusion, let us assume that for a network, the partial graph
$  \mathcal{G}_{0}  $ has a circuit. We may consider the network depicted in
Fig.~\ref{f-FJQN} as an appropriate illustration if we put
$  r_{2}=r_{3}=r_{4}=0 $. In that case, there is the circuit in the graph
$  \mathcal{G}_{0} $, including nodes $  2 $, $  3 $, and $  4 $. Then, as it
follows from Lemma~\ref{l-1}, the implicit dynamic equation associated with
the network cannot be solved in an explicit form.

One can see, however, that it is easy to make the equation solvable only by
setting new initial conditions, without changing the network topology. Since
$  \mathcal{G}_{0} = ({\bf N},{\bf A}_{0}) $, where
$  {\bf A}_{0} = \{(i,j) | \, i \in {\bf P}(j), \, r_{j}=0 \} $, we may
eliminate arcs from the graph $  \mathcal{G}_{0}  $ by substituting nonzero
values for some parameters $  r_{i}  $ set initially to $  0 $. One can
compare the network in Fig.~\ref{f-FJQN}, under the conditions
$  r_{2}=r_{3}=r_{4}=0 $, with that subject to $  r_{2}=r_{4}=0 $, and
$  r_{3}=1 $, as an example.

\bibliographystyle{utphys}

\bibliography{Max-plus_algebra_models_of_queueing_networks}

\end{document}